\documentclass[12pt]{article}

\usepackage{pstricks}
\usepackage{amsmath,amsthm}
\usepackage{amssymb}
\usepackage{color}
\usepackage{xspace}
\usepackage[colorlinks=true,
linkcolor=green,
filecolor=brown,
citecolor=green]{hyperref}
 
\def\gl{ground level\xspace}

\begin{document}
 
\begin{center}
{\Large
Bijections for Dyck paths with all peak heights of the same parity                           \\ 
}

\vspace{2mm}
David Callan   \\
 
\vspace{2mm}

February 17, 2017
\end{center}

\begin{abstract}
We show bijectively that Dyck paths with all peaks at odd height are counted by the Motzkin numbers and Dyck paths with all peaks at even height are counted by the Riordan numbers.
\end{abstract}
\vspace*{2mm}

\section{Introduction} \label{intro}

\vspace*{-4mm}
In a comment on sequence \htmladdnormallink{A005043}{http://oeis.org/A005043} in OEIS \cite{oeis}, Emeric Deutsch notes that Dyck paths of semilength $n$ with all peaks at even height are counted by the Riordan number $R_n$. It is also true that Dyck paths of semilength $n$ with all peaks at odd height are counted by the Motzkin number $M_{n-1}$ (first column of 
\htmladdnormallink{A091869}{http://oeis.org/A091869}).

In this note we give a recursive definition and an explicit description of simple bijections to establish these facts. Of course, due to the well-known glove (or worm-crawls-around-tree) bijection \cite[Fig. 1.8]{cat2015} between ordered trees and Dyck paths, the results can equally well be phrased in terms of ordered trees with all leaf heights of the same parity.
The Motzkin number $M_{n}$ counts Motzkin $n$-paths, that is, lattice paths of $n$ steps formed from upsteps $U=(1,1)$, flatsteps $F=(1,0)$, and downsteps $D=(1,-1)$ such that the path ends at the same level as it starts and never goes below this level (ground level). A Motzkin $2n$-path with no flatsteps is a Dyck $n$-path. Motzkin $n$-paths with no flatsteps at ground level are known to be counted by the Riordan number $R_n$ (\htmladdnormallink{A005043}{http://oeis.org/A005043}).

In Section 2, we present two recursively defined bijections: $\phi_A$ sends 
Dyck $n$-paths with all peaks at odd height to 
Motzkin $n$-paths that start with a flatstep, and $\phi_B$ sends Dyck 
$n$-paths with all peaks at even height to 
Motzkin $n$-paths with no flatsteps at ground level, thereby establishing the stated results. In Section 3, we ``solve the recursion'' and give explicit descriptions of  $\phi_A$ and $\phi_B$; the identical same wording applies for both, only the set they act on differs. 

\section{The bijections defined recursively} \label{recur}
A nonempty Dyck path $P$ with $k$ downsteps returning to \gl can be decomposed as 
\[
UP_1 DU P_2 D \dots UP_k D,
\]
where $P_1,P_2,\dots,P_k$ are all Dyck paths, possibly empty. The paths $UP_i D$ are called the components of $P$. Recall that $\phi_A$ (resp.\ $\phi_B$) acts on Dyck paths all of whose peaks are at odd (resp.\ even) height. First, $\phi_A$ and $\phi_B$ both respect components: 
\[
 \phi_A(P) \textrm{ is the concatenation }\phi_A(UP_1 D)\,\phi_A(U P_2 D) \,\dots \,\phi_A(U P_k D)
\]
and similarly for $\phi_B$. So it suffices to define them on 1-component Dyck paths. Note that if the 1-component Dyck path $UPD$ has all its peaks at odd height, then $P$ has all its peaks at even height and vice versa. The empty Dyck path, denoted $\epsilon$, is considered to have one peak at height 0 (this becomes reasonable under the glove identification of $\epsilon$ with the no-edge one-vertex tree whose root is also a leaf). Now the definitions on 1-component Dyck paths are as follows:
\begin{eqnarray*}
\phi_B(\epsilon)& =&\epsilon\, , \\
\phi_A(U P D) & =& F\, \phi_B(P)\, , \\
\phi_B(U P D) & =& U\, \textrm{Rest} (\phi_A(P))\,D\, ,
\end{eqnarray*}
where $\textrm{Rest} (\phi_A(P))$ denotes $\phi_A(P)$ with its initial step, necessarily $F$, deleted. Clearly, the maps $\phi_A$ and $\phi_B$ are defined on their respective domains and preserve size (number of steps).

The components of a Dyck path $P$ with odd peak heights can be retrieved from 
the Motzkin path $\phi_A(P)$ by picking out the segments that begin with a flatstep at \gl. 
Similarly, the components of a Dyck path $P$ with even peak heights can be retrieved from 
$\phi_B(P)$ by picking out the segments that end with a downstep to \gl.
 
One can easily show by induction on length that the maps are reversible, with respective inverses $\psi_A,\psi_B$ determined by the previous paragraph and the recursion
\begin{eqnarray*}
\psi_B(\epsilon)& =&\epsilon\, , \\
\psi_A(FM) & = & U\, \psi_B(M)\,D\, , \\
\psi_B(U M D) & =& U\, \psi_A(FM))\,D\, ,
\end{eqnarray*}
where $M$ is an arbitrary Motzkin path. 
\section{The bijections described explicitly} \label{explicit}
To present $\phi_A$ and $\phi_B$ explicitly, it is convenient to work with ordered trees all of whose leaf heights are odd (resp.\ even).
Each edge $E$ in such a tree has a \emph{parity}, defined as the parity of the number of edges  (possibly 0) on a path from $E$ (not including $E$) away from the root to a leaf vertex.

\begin{center}
\begin{pspicture}(-10,-1)(10,4)
\psset{unit=.8cm}

\psdots[linewidth=.8pt,linecolor=black](-6,0)
\psdots[linewidth=.5pt,linecolor=black](-11.5,5)(-11,4)(-10.5,5)(-10,3)(-9,4)(-9,5)
\psdots[linewidth=.5pt,linecolor=black](-8,2)(-8,3)(-7,3)(-7,4)(-7,5)(-6,1)(-6,2)(-6,3)(-5,3)(-4,2)(-3,3)(-3,4)(-3,5)
\psline(-6,0)(-6,1)
\psline(-11,4)(-10.5,5)
\psline(-8,3)(-8,2)(-7,3)
\psline(-4,2)(-3,3)

\psset{linewidth=1.0pt}

\psset{linecolor=blue}
\psline(-11,4)(-10,3)(-9,4)
\psline(-8,2)(-6,1)
\psline(-6,2)(-6,1) 
\psline(-4,2)(-6,1)
\psline(-7,4)(-7,3)
\psline(-3,4)(-3,3)

\psset{linecolor=red}
\psline(-11.5,5)(-11,4)
\psline(-8,2)(-10,3)
\psline(-6,2)(-6,3)
\psline(-4,2)(-5,3)
\psline(-7,4)(-7,5)
\psline(-3,4)(-3,5)
\psline(-9,4)(-9,5)

\rput(-6.2,0.5){\textrm{{\scriptsize $a$}}}
\rput(-7.1,1.3){\textrm{{\scriptsize $b$}}}
\rput(-9.1,2.3){\textrm{{\scriptsize $c$}}}
\rput(-10.8,3.4){\textrm{{\scriptsize $d$}}}
\rput(-11.5,4.5){\textrm{{\scriptsize $e$}}}
\rput(-10.5,4.5){\textrm{{\scriptsize $f$}}}
\rput(-9.7,3.6){\textrm{{\scriptsize $g$}}}
\rput(-9.3,4.5){\textrm{{\scriptsize $h$}}}
\rput(-8.2,2.6){\textrm{{\scriptsize $i$}}}
\rput(-7.6,2.7){\textrm{{\scriptsize $j$}}}
\rput(-7.3,3.5){\textrm{{\scriptsize $k$}}}
\rput(-7.3,4.5){\textrm{{\scriptsize $\ell$}}}
\rput(-6.3,1.5){\textrm{{\scriptsize $m$}}}
\rput(-6.3,2.5){\textrm{{\scriptsize $n$}}}
\rput(-5.1,1.7){\textrm{{\scriptsize $o$}}}
\rput(-4.7,2.4){\textrm{{\scriptsize $p$}}}
\rput(-3.6,2.7){\textrm{{\scriptsize $q$}}}
\rput(-3.3,3.5){\textrm{{\scriptsize $r$}}}
\rput(-3.3,4.5){\textrm{{\scriptsize $s$}}}

\rput(4.7,0.5){\textrm{{\scriptsize $a$}}}
\rput(3.5,1.2){\textrm{{\scriptsize $b$}}}
\rput(3.2,2.3){\textrm{{\scriptsize $c$}}}
\rput(.4,2.3){\textrm{{\scriptsize $d$}}}
\rput(0,3.5){\textrm{{\scriptsize $e$}}}
\rput(-1.1,3.5){\textrm{{\scriptsize $f$}}}
\rput(1.1,2.6){\textrm{{\scriptsize $g$}}}
\rput(.7,3.5){\textrm{{\scriptsize $h$}}}
\rput(1.8,2.7){\textrm{{\scriptsize $i$}}}
\rput(2.5,2.8){\textrm{{\scriptsize $j$}}}
\rput(2.7,3.5){\textrm{{\scriptsize $k$}}}
\rput(2.7,4.5){\textrm{{\scriptsize $\ell$}}}
\rput(4.7,1.5){\textrm{{\scriptsize $m$}}}
\rput(4.7,2.5){\textrm{{\scriptsize $n$}}}
 
\rput(5.7,1.8){\textrm{{\scriptsize $o$}}}
\rput(7,2.5){\textrm{{\scriptsize $p$}}}
\rput(6.0,2.4){\textrm{{\scriptsize $q$}}}
\rput(5.7,3.5){\textrm{{\scriptsize $r$}}}
\rput(5.7,4.5){\textrm{{\scriptsize $s$}}}

\psset{linecolor=blue}
\psline(-.5,3)(2,2)(1,3)
\psline(3,3)(3,4)
\psline(6,3)(6,4)
\psline(5,1)(5,2)
\psline(2,2)(5,1)(6.5,2)

\psset{linecolor=red}
\psline(-.5,3)(0,4)
\psline(1,3)(1,4)
\psline(3,4)(3,5)
\psline(2,2)(4.5,3)
\psline(5,2)(5,3)
\psline(6,4)(6,5)
\psline(6.5,2)(7,3)

\psset{linecolor=black}
\psline(5,0)(5,1)
\psline(-.5,3)(-1,4)
\psline(2,3)(2,2)(3,3)
\psline(6.5,2)(6,3)
\psdots[linewidth=.8pt,linecolor=black](5,0)

\psdots[linewidth=.5pt,linecolor=black](-.5,3)(2,2)(1,3)(3,3)(3,4)(6,3)(6,4)(5,2)(5,1)(6.5,2)(0,4)(1,3)(1,4)(3,5)(4.5,3)(5,3)(6,5)(7,3)(2,3)(6,3)(-1,4)

\rput(-2,1){$\longrightarrow$}

\rput(-6.5,-.81){\textrm{\small{Figure 1a}}}
\rput(4,-.8){\textrm{\small{Figure 1b}}}

\end{pspicture}
\end{center}

\begin{center}
\begin{pspicture}(-1,-0.3)(18,2.5)
\psset{unit=.8cm}

\psline[linestyle=dotted](-1,1)(18,1)

\psline(-1,1)(0,1)(1,2)(2,3)(3,3)(4,2)(5,3)(6,2)(7,2)(8,2)(9,3)(10,2)(11,1)(12,2)(13,1)(14,2)(15,2)(16,3)(17,2)(18,1)

\rput(-.5,1.2){\textrm{{\scriptsize $a$}}}
\rput(.4,1.7){\textrm{{\scriptsize $b$}}}
\rput(1.4,2.7){\textrm{{\scriptsize $d$}}}
\rput(2.5,3.3){\textrm{{\scriptsize $f$}}}
\rput(3.6,2.7){\textrm{{\scriptsize $e$}}}
\rput(4.4,2.7){\textrm{{\scriptsize $g$}}}
\rput(5.7,2.7){\textrm{{\scriptsize $h$}}}
\rput(6.5,2.3){\textrm{{\scriptsize $i$}}}
\rput(7.5,2.3){\textrm{{\scriptsize $j$}}}
\rput(8.4,2.7){\textrm{{\scriptsize $k$}}}
\rput(9.6,2.7){\textrm{{\scriptsize $\ell$}}}
\rput(10.6,1.7){\textrm{{\scriptsize $c$}}}
\rput(11.3,1.7){\textrm{{\scriptsize $m$}}}
\rput(12.6,1.7){\textrm{{\scriptsize $n$}}}
\rput(13.4,1.7){\textrm{{\scriptsize $o$}}}
\rput(14.5,2.3){\textrm{{\scriptsize $q$}}}
\rput(15.3,2.7){\textrm{{\scriptsize $r$}}}
\rput(16.6,2.7){\textrm{{\scriptsize $s$}}}
\rput(17.7,1.7){\textrm{{\scriptsize $p$}}}

\rput(-1.5,2.5){$\longrightarrow$}
 
\psdots[linewidth=.5pt,linecolor=black](-1,1)(0,1)(1,2)(2,3)(3,3)(4,2)(5,3)(6,2)(7,2)(8,2)(9,3)(10,2)(11,1)(12,2)(13,1)(14,2)(15,2)(16,3)(17,2)(18,1)
\rput(9,-.1){\textrm{\small{Figure 1c}}}

\end{pspicture}
\end{center}

Given an ordered tree with all peak heights of the same parity, we construct a Motzkin path from it by color coding the edges, tweaking their position in the tree, and then turning the edges into steps according to the color code. The procedure is as follows, with edges/steps labelled in Figure 1 above merely to keep track of them. Color blue each edge of odd parity. Color red the leftmost child edge (there will be one) of each blue edge. Color black all other edges (Fig. 1a). Delete each red edge (in the graph-theoretic sense) and reinsert it as the rightmost child edge of its original parent edge (Fig. 1b). Then walk around the tree clockwise starting at the root (depth first transversal), turning each edge, when first encountered, into a step according to the color code: blue $\rightarrow\, U$, red $\rightarrow\, D$, black $\rightarrow\, F$ (Fig. 1c). 

The example illustrated has odd peak heights but the procedure works just as well for even peak heights. We leave the reader to verify that it coincides with $\phi_A$ for odd peak heights and with $\phi_B$ for even peak heights, and to figure out the inverse of the procedure.

From the definitions, $\phi_A$ sends the statistic ``number of downsteps incident with \gl'' in the Dyck path ( = number of children of root in the corresponding tree) to ``number of \emph{flatsteps} at \gl'' in the Motzkin path, and $\phi_B$ sends it to number of \emph{downsteps} incident with \gl in the Motzkin path.
From the explicit description, one can show that both $\phi_A$ and $\phi_B$ send \#\,peaks in the Dyck path to ($\#\,U$s $-\,\#UU$s)  + ($\#\,F$s  $-\,\#FU$s) in the Motzkin path.

\textbf{Added in Proof.} As Jordan Tirrell 
pointed out to me, the bijections $\phi_A$ and 
$\phi_B$ actually have much simpler descriptions: For $\phi_A$, delete the 
first and last steps of the 
Dyck $n$-path and split the remaining $2n-2$ steps into $n-1$ contiguous pairs. No $UD$ will be 
present when all peaks are at odd height. Then replace the pairs $UU$, $DU$, $DD$ by $U$, $F$, $D
$ respectively and prepend $F$ to get the desired Motzkin path. Similarly for $\phi_B$,  split 
all $2n$ steps into $n$ contiguous pairs. Again, no $UD$ will be present when all peaks are at 
even height. Make the same replacement of pairs to get the desired Riordan path (there will be no 
flatstep at \gl in the Motzkin path because that would require a $DU$ below 
\gl in the Dyck path)\,\cite{tirrell}. 

So this paper is an example of doing things the hard way when there is in fact a much easier way.

\noindent Department of Statistics, University of Wisconsin-Madison, Madison, WI 53706 \ U.S.A. 

\end{document}